\documentclass{amsart}
\usepackage[utf8]{inputenc}
\usepackage [utf8]{inputenc}
\usepackage[math]{iwona}
\usepackage{amsmath,amsfonts, amssymb, version,  enumerate,color,draftcopy,hyperref,graphicx,marginnote,xspace}
\usepackage{layout}

\newtheorem{Theorem}{Theorem}[section]
\newcommand\ol[1]{\overline{#1}}
\newcommand\ip[2]{\langle #1, #2\rangle}
\newcommand\Aut[1]{\aut(#1)}
\DeclareMathOperator{\aut}{Aut}
\newcommand\JBT{JB*-triple\xspace}
\newcommand\tp[3]{\{#1,#2,#3\}}
\newcommand\D[2]{#1 {\mbox{\scriptsize $\square$}} #2}
\newcommand\q[1]{Q_{#1}}
\newcommand\C{\mathbb{C}}
\newcommand\R{\mathbb{R}}
\newcommand\N{\mathbb{N}}
\newcommand\Z{\mathbb{Z}}
\renewcommand\L{\mathcal{L}}
\newcommand\half{\frac12}

\newcommand\inv{^{-1}}

\newcommand\df[1]{{\emph{#1}}\index{#1}}
\newcommand\norm[1]{\lVert #1 \rVert}
\newcommand\summ[3]{\displaystyle\sum_{#1 = #2}^{#3}}
\newcommand\eps{\varepsilon}
	
\newcommand\absval[1]{\vert #1 \vert}

\newtheorem{Lemma}[Theorem]{Lemma}
\newtheorem{Proposition}[Theorem]{Proposition}
\newtheorem{Corollary}[Theorem]{Corollary}

\newtheorem{Def}[Theorem]{Definition}
\newtheorem{Rem}[Theorem]{Remark}
\newtheorem{Ex}[Theorem]{Example}

\newenvironment{Example}{\begin{Ex} \rm}{\end{Ex}\goodbreak}
\newenvironment{Proof}{\noindent\textbf{Proof. }}{\null\hfill$\square$\medskip}
   \newenvironment{Remark}{\begin{Rem}\rm}{\end{Rem}}
\newenvironment{Definition}{\begin{Def}\rm}{\end{Def}\goodbreak}

\DeclareSymbolFont{extrasymbols}{OMS}{cmsy}{m}{n}
\DeclareMathDelimiter{\lVert}
  {\mathopen}{extrasymbols}{"6B}{largesymbols}{"0D}
\DeclareMathDelimiter{\rVert}
  {\mathclose}{extrasymbols}{"6B}{largesymbols}{"0D}

\newcommand\ismc\cong

\title{Holomorphic mappings and their fixed points on Spin Factors}
\author{Michael Mackey}
\email{mackey@maths.ucd.ie}
\address{University College Dublin}
\author{Pauline Mellon }
\email{pauline.mellon@ucd.ie}
\address{University College Dublin}

\newcommand\wto{\overset{w}\to}
\newcommand\wlim{\textrm{w-lim}\ }
\renewcommand\j[1]{j#1}

\begin{document}

\begin{abstract}
In this paper we study holomorphic properties of infinite dimensional spin factors.

Among the infinite dimensional Banach spaces with homogeneous open unit balls, we show that the spin factors are natural outlier spaces in which to ask the question (as was proved in the early 1970s for Hilbert spaces): Do biholomorphic automorphisms $g$ of 
the open unit ball $B$ have fixed points in $\ol B$? In this paper, for infinite dimensional spin factors, we provide reasonable conditions on $g$ that allow us to explicitly construct fixed points of $g$ lying on $\partial B$.  

En route, we also prove that every spin factor has the density property.

In another direction, we focus on (compact) holomorphic maps $f:B\rightarrow B$, having no fixed point in $B$ and examine the sequence of iterates $(f^n)$. As $(f^n)$ does not generally converge, we instead trace the target set $T(f)$ of $f$, that is, the images of all accumulation points of $(f^n)_n$, for any topology finer than
the topology of pointwise convergence on B.  
We prove  for a spin factor that $T(f)$ lies on the boundary of a single bidisc unique to $f$. 
\end{abstract}


\maketitle

\section{Introduction}

Any biholomorphic map on the unit ball $B$ of a Hilbert space extends to
the closed ball $\ol B$ and has a fixed
point there \cite{MR0305158}.   This result by Hayden and Suffridge contrasts with the fact that a homeomorphism of
 $\ol B$ may not have a fixed point \cite{MR0014203}.  It is worthwhile asking which other infinite dimensional complex Banach spaces have this property.   Some Banach spaces have "too few" biholomorphic maps on their open unit ball for this question to be interesting, in particular if every such biholomorphic map is linear and hence, trivially, has a fixed point at the origin.  We focus therefore on spaces with a sufficiently rich group of biholomorphic mappings (or automorphisms) of their open unit ball $B$, by which we mean that the group acts transitively, that is, $\Aut B(0)=B$. These Banach spaces are precisely the \JBT{}s.

 Since it is relatively easy to see that (biholomorphic) automorphisms of a \JBT ball $B$ extend to a neighbourhood of the boundary, our question  becomes: for which \JBT{}s is it true that automorphisms of $B$ must have a fixed point in $\ol B$.
The example in \cite[Counterexample 1.13]{MR1997703}, attributed to Stach\'{o}, shows that the property does not hold for all \JBT{}s as it does not hold in general  for $C(K)$ spaces.  

The above-mentioned result for Hilbert
spaces was essentially proved by appeal to the Brouwer-Schauder fixed point theorem in terms of the weak topology, as direct calculation of the automorphisms showed them  to be weakly continuous and this combined with the weak
compactness of the closed Hilbert ball proved that each automorphism has a fixed point in $\overline B$. 

In general, a \JBT $Z$ does not have a weakly compact ball and indeed, $\ol B$ is weakly compact if, and only if, $Z$  is reflexive if, and only if, $Z$ is finite rank \cite{Kaup_I, Kaup_II}.  The finite rank triples are therefore a natural setting for this question and, to pursue a Brouwer-Schauder approach as before, one might ask whether or not all biholomorphic automorphisms on a finite rank \JBT are weakly continuous.  However, \cite[Corollary 3.8]{MR1166526} shows this is not true generally and, in particular, that on spin factors the only weakly continuous automorphisms of $B$ are the linear ones. 

Still with a view to a Brouwer-Schauder approach, the finite rank triples, being reflexive, are also dual Banach spaces and therefore  have a weak$^*$ topology for which $\ol B$ is  weak$^*$ compact. Indeed, more generally, one could try to use any locally convex topology $\tau$ coarser than the norm toplogy. In \cite[Corollary 3.7]{MR1166526}, it is shown that an automorphism $g$ is $\tau-\tau$ continuous if, and only if, $a=g(0)$ is a $\tau$-continuous element in the sense that the map $:x \mapsto \{x,a,x\}$ is $\tau$ continuous, where $\tau$ is any of the weak (w), weak$^*$ (w$^*$), strong$^*$ or Mackey topologies. A characterisation of all the $\tau$ continuous elements, Cont$_\tau(Z)$, therefore determines whether or not the Hilbert space approach can be used for the $\tau$ topology. It has been shown that Cont$_\tau(Z)$ is a norm closed ideal and it has been explicitly determined for several \JBT s for both the weak and weak$^*$ topologies. However, in almost all of these cases,  Cont$_w(Z)\neq Z$ \cite{MR1073432,MR1166526}. 

In fact, a full characterisation of those triples with  Cont$_w(Z)=Z$ is given in \cite[Theorem 5.7]{MR1273340}.  In particular, for finite rank triples $Z$ satisfying the complicated three factor condition given in \cite[Theorem 5.7]{MR1273340} (namely $Z$ has dual RNP, all of its spin factor representations are finite dimensional and every w$^*$-dense representation onto a Cartan factor is elementary), the Browder-Schauder approach works and proves  that all automorphisms of such $Z$ do indeed have fixed points in $\ol B.$ It also shows, however, that a similar approach is not possible in general.

Spin factors are natural outliers here; in particular, if $J$ is a spin factor then Cont$_w(J)=\mbox{Cont}_{w^*}(J)=\{0\}$ \cite[Proposition 4.3]{MR1073432} and\cite[Corollary 3.7]{MR1166526}. It makes sense therefore to focus on the spin factors as being somehow "least likely" candidates for a Hayden-Suffridge type fixed point result, necessitating also an entirely different approach.

This paper therefore starts by explicitly calculating a biholomorphic automorphism $g$ on the unit ball $B$ of a Spin Factor $J$. The weak compactness of $\ol B$ is then used to obtain a point $\xi \in \ol B$ that we call a  {\it{weak fixed point}} of $g$. While this weak fixed point may or may not be an actual fixed point, our first main result Theorem~\ref{thm:first} gives explicit conditions on the weak fixed point $\xi$ that allow us to conclude that the weak fixed point is indeed a fixed point.  Thereafter, we focus on the situation where these conditions fail and proceed by constructing a fixed  point by different means.   In general, the constructed fixed point may not be equal to the necessarily existing weak fixed point $\xi$.    Subsequent results, Theorems \ref{thm:orthogonal}
and \ref{thm:sliver_condition},  then reduce the conditions required for $g$ to have a fixed point. Unlike the Hilbert space proof which relies on existence results, our approach is constructive, namely, we solve for the fixed point, which we show belongs to $\partial B$.

The final short section of the paper returns to earlier work of the authors on iterates of compact holomorphic maps $f:B \rightarrow B$ which have no fixed point in $B$. As it is known that the sequence of iterates  $(f^n)_n$ does not generally converge, a concept of Wolff Hull, $W(f)$, was introduced in \cite{MackeyMellon_WH} to help locate the images of accumulation points $\Gamma(f)$ of $(f^n)_n$, for any topology finer than the topology of point-wise convergence on $B$, that is, to locate the target set $$T(f):=\bigcup_{g \in \Gamma(f)} g(B).$$

For a spin factor $J$ and $f:B \rightarrow B$ as above, Theorem~\ref{Wolff hull} and Corollary~\ref{bidisc} below show that there exists $e \in \partial B$ and unique bidisc $B_e:=\Delta e \times \Delta \j e$ (where $j$ is the conjugation on $J$) such that  $$T(f) \subseteq \partial B_e.$$  We note that $B_e$ is the open unit ball of $\C e \times \C \j e$, which is a copy of $\C^2$ with the maximum norm. This result is best possible in light of \cite[Theorem 3.4]{Mellon-dynamics}.

 Another small but novel result included here is Corollary~\ref{density} below which shows that spin factors have the density property (cf. \cite{Densprop}).

\section{\JBT{}s and Spin Factors}
\label{sec:JBT}

If the open unit ball $B$ of a complex Banach space $Z$ has a transitive group of biholomorphic mappings acting on it, i.e.~$\Aut B(0)=B$, then the symmetry at the origin, $s_0(x)=-x$, provides a symmetry at every point via conjugation: $s_a= g_a\circ s_0\circ g_a\inv$, where $g_a$ is any automorphism taking the origin to $a$.  A \df{bounded symmetric domain} is any domain in a complex Banach space which is biholomorphically equivalent to such an open unit ball and then, and only then, does the Banach space $Z$ possess the algebraic structure that renders it a \JBT.  

\begin{Definition}\label{def:JBT}
    A \JBT\ is a complex Banach space $Z$ with a real trilinear
    mapping $\tp\cdot\cdot\cdot:Z\times Z\times Z \to Z$ satisfying
    \begin{enumerate}[(i)]
        \item $\tp xyz$ is complex linear and symmetric in the
        outer variables $x$ and $z$, and is complex anti-linear
        in $y$.

                \item The map $z\mapsto \tp xxz$, denoted $\D xx$, is
                Hermitian, $\sigma(\D xx)\ge 0$ and
        $\norm{\D xx} = \norm x^2$ for all $x\in Z$, where
        $\sigma$ denotes the spectrum.

        \item The product satisfies the following ``Jordan triple
        identity''
            \[\tp ab{\tp xyz} = \tp{\tp abx}yz -
              \tp x{\tp bay}z + \tp xy{\tp abz}.\]
    \end{enumerate}
\end{Definition}
 
The triple product is jointly continuous 
($\norm{\{x,y,z\}} \leq \norm{x} \norm{y} \norm{z}$ \cite{FR_GN}) and gives rise to the bounded linear maps: 
 $\D xy \in \L(Z):z\mapsto \tp xyz$,
    $\q x \in \L_{\R}(Z):z\mapsto \tp xzx$,
and the  Bergman operators   $B(x,y)=I-2\D xy + \q x\q y \in \L(Z)$.  For $a\in B$, the Bergmann  operator $B(a,a)$ has positive spectrum and
a unique square root  $B(a,a)^\half$ with positive spectrum, denoted by $B_a$, and this satisifies $\norm {B_a\inv}=\frac 1{1-\norm a^2}$ and $\norm{B_a}\le 1$ (see \cite{Kaup_Oviedo} or \cite{MackeyMellon_CO}).  A key fact is that $T\in \L(Z)$ is a surjective isometry if, and only if, it is a bijective homomorphism of the triple product.  

   Definition~\ref{def:JBT}  (cf. \cite{Kaup_RMT}) allows 
an explicit description of the biholomorphic automorphisms of $B$, denoted $\Aut B$, namely,
every $g\in\Aut B$ can be written uniquely in the form \begin{equation}\label{eq:automorphism}
    g=Tg_a,
\end{equation} where
$T$ is a surjective linear isometry of $Z$ and $g_a$ is a generalised
Möbius map, or \df{transvection}, defined on $B$ by
\begin{equation}\label{eq:mobius}
    g_a(x)= a+ B_a x^{-a}
\end{equation}
where $x^{a}:=(I-\D xa)^{-1}x$ is known as the quasi-inverse of $x$ with respect to $a$.  Indeed $(x,y)\in Z^2$ is said to be a quasi-invertible pair if $B(x,y)$ is invertible and then $x^y=B(x,y)\inv(x-\q xy)$ is the quasi-inverse of $x$ with respect to $y$\ \cite{Densprop}. 
Evidently, $g_a(0)=a$ and  $g_a\inv= g_{-a}$.  For
$a\in B$, the quasi-inverse map $x\mapsto x^{-a}$, and hence $g_a$ and every
element of $\Aut B$, is defined and
continuous beyond the unit ball, to the ball of radius $\norm a\inv$.

A non-zero element $e \in Z$ is a \df{tripotent} if $\{e,e,e\}=e$. Every tripotent $e$ induces a splitting of $Z$, as
	$Z=Z_0(e) \oplus Z_{\frac{1}{2}}(e) \oplus Z_1(e)$,
where $Z_k(e)$ is the $k$ eigenspace of $e \square e$ and the linear maps 
$P_0(e)=B(e,e), P_{\frac{1}{2}}(e)=2(\D ee -\q e\q e),$ and $ P_1(e)=\q e\q e$ are mutually orthogonal projections of $Z$ onto $Z_0(e),\  Z_{\frac{1}{2}}(e),$ and $Z_1(e)$ respectively.
 
Elements $x,y \in Z$ are (triple) orthogonal, $x \bot \ y$, if $x \square y = 0$
(or equivalently \cite{loos_bsd} if $y \square x = 0$). A tripotent is \df{minimal} if $Z_1(e)=\C e$, \df{maximal} if $Z_0(e)=\{0\}$ and \df{unitary} if $Z_1(e)=Z$.  If $c$ and $e$ are orthogonal tripotents then $c+e$ is also a tripotent.  Two tripotents $e$ and $f$ will be orthogonal if $\tp eef =0$.  If a tripotent $e$ is the sum of $r$ mutually orthogonal minimal tripotents then it is said to have rank $r$ and this is well-defined.  

Examples of \JBT{}s include the space of bounded linear operators between Hilbert spaces, $\L(H,K)$, where \begin{equation}\tp xyz=\half(xy^*z+zy^*x),\label{eq:jtp}\end{equation} and their norm closed subtriples.  In particular, all C*-algebras are JB*-triples with respect to this triple product, as is a Hilbert space $H$ via identification with $\L(\C, H)$ and the resulting triple product $\tp xyz=\half(\ip xy z+ \ip zy x)$. 
Throughout the rest of this paper we concentrate on the \JBT{s} known as \df{spin factors}.

\section{Biholomorphic Automorphisms of Spin Factors}
\label{sec:sf}

There are different, but equivalent, ways of introducing spin factors  (see \cite{MR1195079} for a grid approach and \cite{Din_Sch} for a matrix treatment), which in the literature may also be termed \df{spaces of spinors} or \df{Cartan factors of Type IV}.   Our approach is quite standard, which is to define a spin factor $J$ as a Hilbert space equipped with
conjugation $j:J\to J$ (that is, $j$ is conjugate linear, $\ip {\j x}{\j y}=\ol{\ip xy}$ and $j^2=I$) under the equivalent norm
\[ \norm x^2 = \ip xx +\sqrt{\ip xx^2-\absval{\ip x{\j x}}^2}\]
and triple product
\begin{equation}
\tp abc = \ip ab c+ \ip cba -\ip
a{\j c}\j b. \label{eq:tp}
\end{equation}   Note that, although the conjugation is not linear, it is both isometric, $\norm{\j x}=\norm x$, and a triple product preserver, $j\tp abc = \tp{ja}{jb}{jc}$.  

By definition a spin factor is isomorphic to a Hilbert space and consequently is reflexive and falls into the category of JBW*-triples, which are the \JBT{}s that have a Banach space predual.   Spin factors enjoy the Kadec-Klee property \cite[Corollary 3]{MR1874485} which is to say that if $x_n\to x$ weakly and $\norm{x_n}\to \norm{x}$ then $x_n\to x$.

\subsection{Tripotents}
A tripotent in a spin factor satisfies \[ 0\ne e=\tp eee = 2\ip eee -\ip e{je}je\]
so either
\begin{enumerate}
\item[(i)] $\ip e{je}=0$ and $\ip ee = \half$ (and $e$ is minimal) or
\item[(ii)] $je=\lambda e$ with $|\lambda|=1$ and $\ip ee=1=\lambda\ip e{je}$ (and $e$ is maximal).
\end{enumerate}
Indeed in (i), $\tp exe=2\ip ex e$ so $P_1^e(x)=\q e\q e x =2\ip x
e e \in \C.e$ and $e$ is a minimal (rank 1) tripotent.  While in  (ii), $\q
e\q e x=x$ and $e$ is a maximal tripotent which has rank 2 (see below).  Maximal
tripotents in a spin factor are unitary.

We note that tripotent orthogonality is distinct from
inner product orthogonality, as follows.

\begin{Lemma}\label{lem:orthmintp}
  Let $e$ and $f$ be minimal tripotents in a spin factor.  Then $e$ and $f$ are triple orthogonal if, and only if,
  \begin{enumerate}[(a)]
  \item $\ip ef=0$, and
    \item $jf=\lambda e$ for
      some $|\lambda|=1$ (and then $je=\lambda f$).
  \end{enumerate}
\end{Lemma}

\begin{Proof}
  First, if $e\perp f$ then $0=\tp efe =2\ip ef e - \ip e{je}f = 2\ip ef e$ and  so $\ip
  ef =0$. Also we have, \begin{align*}
                                      0=\tp eff &=\ip eff + \ip ffe
                                                 -\ip e{jf}jf\\
                                      &= \half e - \ip e{jf}jf
                                    \end{align*}
    which shows $jf= \lambda e$.  Taking norms shows $|\lambda|=1$.

    Conversely, (a) and (b) together show that $\tp eef = \ip eef + \ip fee - \ip e{jf}je= f/2 + 0 -\ol \lambda \ip ee\lambda f =0 $ and hence $e\perp f$.
  \end{Proof}

\begin{Remark}\label{remark:basis}
Given the Hilbert space $H$ and a conjugation $j$, one may choose an orthonormal basis of $H$ for which $j:\sum
\alpha_ke_k \mapsto \sum \ol{\alpha}_k e_k$.    The elements of this orthonormal basis are tripotents but are not \emph{orthogonal tripotents} of the spin factor $J=(H, j)$.  In fact, each 
$e_k$ is a maximal tripotent fixed by $j$ and for any pair of indices $k\ne m$, the elements $(e_k+ie_m)/2$ and
$(e_k-ie_m)/2=j(\frac{e_k+ie_m}2)$ are orthogonal minimal tripotents whose sum equals $e_k$.  Conversely, if $d$ is a minimal tripotent then so is $jd$, and is orthogonal to $d$ while $d+jd$ and $d-jd$ are maximal tripotents.  
\end{Remark}

\begin{Lemma}\label{lem:maxtp}
  For $x\in J$, $jx=\alpha x$ for $\alpha\in \C$ (necessarily, $|\alpha|=1)$ if, and only if,
  $x$ is a scalar multiple of a maximal tripotent.
\end{Lemma}

\begin{Proof}
  If $x=\beta e$ for a maximal tripotent $e$ then $jx=\ol\beta je$.
  As $e$ is a maximal tripotent $je=\lambda e$, where $|\lambda|=1$
  and so $jx =(\frac{\ol \beta}{\beta} \lambda) \beta e = \mu x$ where
  $|\mu|=1$.

  Conversely, if $jx=\alpha x$ then necessarily $|\alpha|=1$ by
  consideration of norms.  It follows $\tp xxx = \ip xxx$ so that
  $y=\frac x{\sqrt{\ip xx}}$ is a tripotent.  This tripotent $y$ is
  maximal since $\ip yy =1$.    
\end{Proof}

\begin{Lemma}\label{lem:mintp}
  For $x\in J$, $\ip x{jx}=0$ if, and only if, $x$ is a scalar multiple of a minimal tripotent.
\end{Lemma}

\begin{Proof}
  If $0\ne x=\lambda e$ for a minimal tripotent $e$ then
  $\ip e{je}=0$, and so too $\ip x{jx}=0$.  Conversely, if $\ip
  x{jx}=0$ then $\tp xxx = 2\ip xx x$ and $\frac x{\sqrt{2\ip xx}}$ is
  a minimal tripotent.
\end{Proof}

Remark that a scalar multiple of a tripotent is a non-negative multiple of a (different) tripotent of the same rank because the unimodular part can be taken into the tripotent, i.e. $\lambda e= r\exp(i\theta)e = r e'$ where $e'=\exp(i\theta)e$ is also a tripotent that is maximal precisely when $e$ is.

\begin{Corollary}\label{cor:3.4}
  $\norm x^2= \ip xx$ (i.e., $\norm x_J = \norm x_H$) if, and only if, $x$ is a
  scalar multiple of a maximal tripotent.
\end{Corollary}

\subsection{Quasi-invertibility in the spin factor}

\begin{Lemma}\label{lem:sfqi}
  For $x,y \in B$ the quasi-inverse $x^y$ is given by \[ x^y= \frac
    1{1-2\ip xy + {\ip x{jx}\ip{jy}y}} \left(x - {\ip x{jx}}jy\right).\]
\end{Lemma}

\begin{Proof}
  Using the triple product of \eqref{eq:tp}, we find that
  \begin{equation}
\tp xyx =
  2\ip xyx-{\ip x{jx}} jy = \gamma x-\xi jy.\label{eq:1}
\end{equation}
where $\gamma=2\ip xy$ and $\xi={\ip x{jx}}$.  Further,
\begin{equation}
  \label{eq:2}
  \tp xy{jy} = {\ip{jy}y} = \eta x
\end{equation}
for $\eta={\ip{jy}y}$. It follows that
\[ (\D xy)(\alpha x+\beta jy) = (\alpha\gamma - \beta\eta)x -\alpha\xi
jy\]
and $\D xy$ is invariant on the span of the basis $\{x, jy\}$,  where
it is represented by
the matrix $T=
\begin{bmatrix}
  \gamma &\eta\\ -\xi & 0
\end{bmatrix}$.  Since $x^y= (I-\D xy)\inv x = \summ k0\infty (\D xy)^k x$, we find
$x^y$ is represented with respect to this basis by $(I-T)\inv
\begin{bmatrix}
  1\\0
\end{bmatrix}$, that is 
\begin{align}
  x^y&= \frac 1{1-\gamma+\xi\eta}
  \begin{bmatrix}
    1 & \eta\\ -\xi &1-\gamma
  \end{bmatrix}\begin{bmatrix}
    1\\0
  \end{bmatrix} \\
&= \frac 1{1-\gamma+\xi\eta} \begin{bmatrix} 1\\-\xi
  \end{bmatrix} \\
&= \frac
    1{1-2\ip xy +  {\ip x{jx}\ip{jy}y}} \bigl(x -  {\ip x{jx}}jy\bigr).
\end{align}
The argument is simplified and the formula remains valid (but simplifies) if $x$ and $jy$ are linearly dependent.
\end{Proof}

\subsection{Bergmann Operators and the Density Property}

\newcommand\ii{\ip x{jx}\ip{jy}y}
\newcommand\linspan{{\rm{span}}}
Recall that $B(x,y)z= z-2\tp xyz + \tp x{\tp yzy}x$.
Calculations similar to those above give
\begin{align} B(x,y)z &= x\left[-2\ip zy(1-2\ip xy) -2 \ii  \right] \notag
  \\
                       &\qquad +jy \left[ 2 \ip x{jz}-2 \ip x{jx}\ip zy \right] \notag\\
  &\qquad\qquad +z \left[1- 2 \ip xy +\ii\right]. \label{eq:berg}
\end{align}
In particular, $B(x,y)z$ lies in the linear span of the elements $\{x,jy,z\}$.  Moreover,
\begin{align}
  B(x,y)x &= x\left( (1- 2\ip xy)^2- \ii \right) +jy\left( 2\ip
            x{jx}(1-\ip xy)\right)\label{eq:12}
            \intertext{and}
            B(x,y)jy &= x\left( -2\ip{jy}y(1-\ip xy) \right) +jy
                       \left( 1-\ii \right).\label{eq:13}
\end{align}

\begin{Proposition}\label{prop:berginv}
  The Bergmann operator $B(x,y)$ is invertible if, and
  only if, $$r(x,y):= 1-2\ip xy + \ii \ne 0.$$
\end{Proposition}

\begin{Proof}
  Suppose $r(x,y)\ne 0$. We show $B(x,y)$ is invertible.  First
  injectivity and suppose, for sake of contradiction, that for some
  $z\ne 0$,  we have $B(x,y)z=0$.  If  $z \notin \linspan\{x, jy\}$  then
  \eqref{eq:berg} provides the contradiction that $r(x,y)=0$.  On the
  other hand, if $z\in \linspan\{x, jy\}$, so $z=\alpha x+ \beta jy$,
  then with respect to the basis $\{x, jy\}$, \eqref{eq:12} and \eqref{eq:13} show that
  \[
    \begin{bmatrix}
      (1-2\ip xy)^2-\ip x{jx}\ip{jy}y  &-2\ip{jy}y(1-\ip xy) \\ 2\ip x{jx}(1-\ip xy) &1-\ip x{jx}\ip{jy}y
    \end{bmatrix}
    \begin{bmatrix}
      \alpha \\ \beta
    \end{bmatrix} =
    \begin{bmatrix}
      0 \\ 0
    \end{bmatrix}.
\]
Thus the determinant of this matrix vanishes, but this determinant is precisely
$r(x,y)^2$ and we have a contradiction.  We conclude that $r(x,y)\ne 0$ implies injectivity.  In
fact, this argument is easily seen to be reversible and so $B(x,y)$ is
injective precisely when $r(x,y)\ne 0$.

Surjectivity of $B(x,y)$ goes along similar lines.  If $r(x,y)\ne
0$ then, given $z\in J$, we can use \eqref{eq:berg}-\eqref{eq:13} and elementary linear algebra (this time we encounter a $3\times 3$ matrix whose determinant is $r(x,y)^3$) to
find $w\in J$, and indeed $w\in\linspan\{x, jy, z\}$ such that $B(x,y)w=z$.  

\end{Proof}

The \df{density property} for a \JBT is said to hold when the set of quasi-invertible pairs $\{(x,y): B(x,y)\mbox{ is invertible}\}$ is dense in $Z\times Z$.  This property has relevance to the construction of a compact type symmetric manifold modelled on $Z$, as described in \cite{loos_bsd}.  While a list of \JBT{s} with the density property was given in \cite{Densprop} (see also \cite{MR1997703}) the spin factors were not considered there. However, the following unpublished result appeared in \cite{Keenan_project}.

\begin{Corollary}\label{density}
  Spin factors have the {density property} \cite{Densprop}, that
  is, for any $x,y \in J$ and $\eps>0$, there exists $z\in J$ with
  $\norm{z}<\eps$ such that $B(x,y+z)$ is invertible.  
\end{Corollary}

\begin{Proof}
    If the result fails then there exist $x,y \in J$ and $\eps>0$ such that for all $z\in B$ and $0\le t<\eps$, the operator $B(x,y+tz)$ is not invertible and, by Proposition~\ref{prop:berginv}, then $r(x,y+tz)=1-2\ip x{y+tz} + \ip x{jx}\ip{j(y+tz)}{y+tz}=0$.  For all $t\in(0,1)$ we then have
    \begin{align*} 0 &= \frac 1t\left(r(x,y+tz)-r(x,y)\right) \\
    &= -2\ip xz + \ip x{jx}(\ip{jy}z + \ip{jz}y + t\ip{jz}z).
    \end{align*}
    This implies that $\ip{jz}z=0$ for all unit vectors $z$, and by Lemma~\ref{lem:mintp}, that every unit vector is a minimal tripotent which is clearly false (by Remark~\ref{remark:basis} for example).
\end{Proof}

\subsection{Linear automorphisms}\label{sec:linaut}
A surjective linear isometry $T$ of $H$ satisfies $\ip{Tx}{Tx}=\ip xx$
but may not be a triple isometry of $J=(H,j)$.  Indeed, a surjective linear isometry of $H$ may not preserve the triple rank of elements in
$J$, while rank must be preserved by a JB*-triple automorphism.  However, it is easy to
check that if $U$ is a surjective linear isometry of $H$ (i.e. a
unitary operator), and $Uj=jU$ then it follows that for the triple norm
$\norm {Ux}=\norm x$.     In fact, up to a scalar, this characterises
triple automorphisms of $J$:  if $T$ is a surjective linear
isometry of a spin factor $J$ then
$T=\exp(i\theta)U$, where $U$ is a unitary operator on $H$ commuting
with $j$, see \cite[p. 196]{MR1195079}.

\begin{Example}
    
Let $(e_n)_n$ be an orthogonal basis  of $H$ with $j(e_n)=e_n,\ \hbox{for all}\ n \in \N.$
\begin{itemize}
\item[(i)] A permutation of the basis, and more generally, a unimodular scalar multiple of basis
  permutation, is a surjective linear isometry of $J$.
  \item[(ii)] The unitary operator $W$ sending $e_1\mapsto \frac 1{\sqrt
      2}(e_1+ie_2)$ and $e_2 \mapsto \frac 1{\sqrt2}(e_1-ie_2)$,
  and otherwise fixing the basis, is a linear isometry of $H$ but not
  of $J$.  It preserves neither the spin factor norm nor the rank of the maximal
  tripotent $e_1$ and clearly $W$ does not commute with $j$ as $jWe_1 = j(e_1+ie_2)/\sqrt 2
  =(e_1-ie_2)/\sqrt 2 \ne We_1 = Wj e_1$.
  \item[(iii)] The even-simpler looking unitary operator $V$ sending
    $e_1\mapsto ie_1$ and $e_k\mapsto e_k$ for $k>1$ is not an
    isometry of $J$, since from
    Lemma~\ref{lem:mintp} above,
    rank$(ie_1+e_2)=1$ while, from Lemma~\ref{lem:maxtp}, rank$(e_1+e_2)=2$.
\end{itemize}
In particular, if a surjective linear isometry $T$ of $J$ commutes with the conjugation $J$ then it is a unitary operator and $\ip{T^kx}{T^ky}=\ip xy$, for all $k\in \Z$.

\end{Example}
\subsection{Möbius maps}

We develop explicit formulae for the M\"obius map, or transvection, $g_a$, given in  (\ref{eq:mobius}) above for $a$  a scalar multiple of a tripotent.  The case of a maximal tripotent is of particular use later.

\begin{Corollary}\label{cor:gmax}
Let $e$ be a maximal tripotent with $je=\sigma e$, $|\sigma|=1$ and let $a=te$, $t\in [0,1]$.
   Then
  $$g_a(z)=\left(t + \frac{\sigma t(1-t^2)\ip z{jz}}{1+2t\ip ze +
      \sigma t^2\ip z{jz}} \right)e + \frac{1-t^2}{1+2t\ip
    ze+\sigma t^2\ip z{jz}} z.$$
\end{Corollary}

\begin{Proof}
  Let $a=te$ as stated. Since $e$ is a maximal tripotent $\ip ee=1$ and
  $\D ee =Id=\q e\q e$ and $B(e,e)=0$.  It follows that  $B(te,te)=(1-t^2)^2Id$
  and so from equation (\ref{eq:mobius}) above 
  $$g_a(z)= a+B(a,a)^\half(z^{-a})=te+(1-t^2)\frac{z+t\ip
    z{jz}je}{1+2t\ip ze +\ip z{jz}t^2\ip {je}{e}}$$
and the stated formula immediately follows. 
\end{Proof}

\begin{Corollary}\label{minaut}
  Let $a$ be a scalar multiple of a minimal tripotent $e$ so that
  $a=te$, $t\in [0,1]$, $\ip  e{je}=0$, $\ip ee = \half$.  Then
  \[g_a(z)=\frac{t+\zeta}{1+t\zeta}e  + \frac{\zeta'+t\ip
      z{jz}}{1+t\zeta} je + \frac{\sqrt{1-t^2}}{1+t\zeta}\tilde z\]
  where $\zeta=2\ip ze$, $\zeta'=2\ip z{je}$ and $z=\zeta e+ \zeta'je
  + \tilde z$ with $\ip{\tilde z}e=\ip{\tilde z}{je}=0$.
\end{Corollary}

\begin{Proof}
  Since $e$ and $je$ are Hilbert orthogonal, we may extend from these two
  basis vectors to an orthogonal basis of $H$.   For $z \in J$, write
  $z=\zeta e
  +\zeta'  je + \tilde z$ where $\ip e{\tilde z}=\ip {je}{\tilde z}=0$ and $\zeta=2\ip ze$,
$\zeta'=2\ip z{je}$.

Standard calculations now show that $B(te,te)e=(1-t^2)^2e$,
$B(te,te)je = je$ and $B(te,te)\tilde z=(1-t^2)z$.  Hence
\[B(te,te)z= (1-t^2)^2\zeta e + \zeta'
je + {(1-t^2)}\tilde z\] and 
\begin{equation}\label{eq:bb} B(te,te)^\half z =  (1-t^2) \zeta e + \zeta'
je + \sqrt{1-t^2}\tilde z.
\end{equation}
From Lemma~\ref{lem:sfqi},
\begin{align*}
z^{-a}=z^{-te}&= \frac 1{1+2t\ip ze+\ip z{jz}t^2\ip{je}e}(z+t\ip
                z{jz}je)\\
  &= \frac 1{1+t\zeta}(\zeta e +(\zeta'+t\ip z{jz})je + \tilde z)
\end{align*}
and, using (\eqref{eq:bb}) we have
\begin{align*}
  g_a(z)&= a+B(a,a)^\half z^{-a} \\
      &=te+ \frac 1{1+t\zeta} B(te,te)^{\half}(\zeta e +(\zeta'+t\ip z{jz})je +
        \tilde z) \\
  &= te+ \frac 1{1+t\zeta}((1-t^2)\zeta e +(\zeta' +t\ip z{jz}) je +
    \sqrt{1-t^2}\tilde z) \\
  &= \frac{t+\zeta}{1+t\zeta}e  + \frac{\zeta' +t\ip z{jz}}{1+t\zeta}
    je + \frac{\sqrt{1-t^2}}{1+t\zeta}\tilde z
\end{align*}
\end{Proof}

\section{Fixed points of Biholomorphic Automorphisms}
\label{sec:fp}

For all $JB^*$-triples each biholomorphic automorphim $g \in \Aut B$ is given, cf. (\ref{eq:automorphism}),  by
$g=Tg_a$, for $T$ a surjective linear isometry  and $g_a$ a Möbius map, for some $a\in B$. While this decomposition is unique, we also have  $g=Tg_a = g_{Ta} T$ \cite{Kaup_RMT}.   In particular,\begin{equation}\label{exp}
\exp(i\theta)g_a = g_{\exp(i\theta) a}\exp(i\theta)\ \hbox{for}\  \theta \in \R.\end{equation} Note first that although both $T$ and $g_a$  separately have fixed points in $\ol B$, given by $0$ and the support tripotent of $a$ respectively \cite[Theorem 0.1]{Mellon-dynamics}, this, of course, does not guarantee that $g$ has a fixed point.  

\begin{Remark}
\begin{enumerate}[(i)]
    \item If $T$ is invariant  the smallest closed subtriple of $Z$ containing $a$ (denoted $Z_a$) which is finite dimensional for any finite rank $Z$, then finite dimensional arguments apply to guarantee the existence of a fixed point of $g$.  In particular, if $T$ is a scalar rotation then $g=Tg_a$ has a fixed point.
    \item Similarly, if $T^n=I$ for some $n$ then the subtriple $V_a$ generated by $\{a,Ta,T^2a,\ldots,T^{n-1}a\}$, being finitely generated, is finite dimensional and is invariant under $g$.  Therefore, again by finite dimensionality, $g$ has a fixed point in $V_a\subset Z$. 

\end{enumerate}

\end{Remark}

We return now to a spin factor $J$ and $g \in \Aut B.$ As described earlier, $Cont_w(B)=\{0\}$, and hence if $g$  is weakly continuous it is necessarily linear. In particular, we cannot use a Brouwer-Schauder approach to discover whether or not a non-linear $g$ has a  fixed point.  We therefore explore an alternative constructive approach. 
    
A common technique to locate fixed points for a holomorphic function on a domain in a Banach space is to combine use of the Earle-Hamilton theorem with compactness (in some topology) of either the domain or the function.   Here, we choose $\alpha_k \uparrow 1$ and let $z_k \in B$ be a fixed point
of $\alpha_kg$ as guaranteed by the Earle-Hamilton theorem.  Since $\ol B$ is weakly sequentially compact, we can assume without loss of generality that
\begin{equation}
z_k \wto \xi \in \ol{B}.\label{eq:0}
\end{equation}

A point $\xi$ in $\ol{B}$ which is a weak limit of a sequence $(z_k)_k$ satisfying $z_k=\alpha_k g(z_k)$ with $\alpha_k\to 1$ will be referred to here as a  {\it{weak fixed point}} of $g\in
\textrm{Aut}(J)$.  
While $g \in \Aut B$, and indeed any holomorphic self-map of $B$, has a weak fixed point, Example~\ref{ex:helpful} below shows that this weak fixed point is not necessarily a fixed point. 

Of course, in a Hilbert space, each $g \in \Aut B$ is weakly continuous and hence
$\xi=\wlim z_k= \wlim \alpha_k g(z_k)=\wlim g(z_k)=g(\wlim z_k)=g(\xi).$

\begin{Proposition}\label{prop:zero}
  Let $J$ be a spin factor, $B$ its open unit ball, $a\in B$ and $(z_k)_k\subset B$.  Suppose $z_k\wto 0$ and $g_a(z_k) \wto 0$, where $g_a$ is the transvection of $B$ generated by $a$. 
  Then $a$ is a scalar multiple of a maximal tripotent. 

    Moreover, either $a=0$ or $\lim \norm{z_k}=\lim \norm{g_a(z_k)} =1$.
\end{Proposition}

\begin{Proof}
  Since $g_a(x)=a+B(a,a)^{\half} x^{-a}$, we have
  \begin{equation}
(z_k)^{-a} \wto
  B(a,a)^{-\half}(-a).\label{eq:wlim}
\end{equation}

In any JB*-triple, with $a$ in the ball, $B(a,a)^{-\half}(a)= a^a$
and, via calculation in the local triple generated by $a$,
$\norm{a^a}=\frac{\norm a}{1-\norm a^2}$.  In $J$ we have
by \ref{lem:sfqi} that \begin{equation} B(a,a)^{-\half}(-a)=
  \frac{-1}{1-2\ip aa + \absval{{\ip a{ja}} }^2} \left( a- {\ip a{ja}}
    ja \right).
\label{eq:one}\end{equation}
Also by \ref{lem:sfqi}, 
\[ (z_k)^{-a}=  \frac
    1{1+2\ip {z_k}a +  {\ip {z_k}{jz_k}\ip{ja}a}} \left(z_k +  {\ip {z_k}{jz_k}}ja\right).
\]  and, choosing a subsequence of the weakly null sequence $z_k$ so that $\ip{z_k}{jz_k}$
converges (to $s$ say), we see that the weak limit of $(z_k)^{-a}$ is
\[ (z_k)^{-a}\wto  \frac
    1{1+  {s\ip{ja}a}} \left( s.ja\right).\]
Uniqueness of limits  together with \eqref{eq:wlim} and
\eqref{eq:one} gives
\begin{equation}
  \label{eq:a}
  \frac{-1}{1-2\ip aa + \absval{{\ip
        a{ja}}{} }^2} \left( a- {\ip a{ja}} ja \right) = \frac
    1{1+  {s\ip{ja}a}} \left( s.ja\right).
\end{equation}
This implies that $a$ and $ja$ are
linearly dependent so $ja=\lambda a$.  Lemma~\ref{lem:maxtp} then implies
$a$ is a scalar multiple of a maximal tripotent.  

Without loss of
generality, we assume  $\lambda=1$.  Indeed, $\norm a =\norm {ja}$ gives $|\lambda|=1$.  Let $\tau^2=\lambda$.  Then from (\ref{exp}) $\tau g_a(z_k)=g_{a'}(z_k')$ where $a'=\tau a$ and $z_k'=\tau z_k$.  Evidently, $j(a')=j(\tau a)=\ol\tau ja =\ol \tau \tau^2 a= \tau a=a'$ and both $z_k'$ and $g_{a'}(z_k')$ converge weakly to $0$. This justifies our assumption, since $a$ is a scalar multiple of a maximal tripotent if, and only if, $a'$ is also.

Remark that \eqref{eq:a} now becomes
\begin{align*}
\frac{-1}{(1-\ip aa)^2} \left( a- \ip a{a} a \right) &= \frac
    1{1+  s \ip{a}a} \left(  sa\right)\\
\frac {-1}{1-\ip aa}a &= \frac{ s}{1+ s\ip aa} a  
\end{align*}
and so either $a=0$ (and we are done) or $s =-1$, that is, $\ip{z_k}{jz_k} \to
-1$, which via the Cauchy Schwarz inequality implies $\norm{z_k}\to 1$.  Indeed, $1\ge \norm{z_k}^2\norm{jz_k}^2 \ge \ip{z_k}{z_k}\ip{jz_k}{jz_k} \ge |\ip{z_k}{jz_k}|^2 \to 1$, and hence $\norm{z_k}\to 1$.  Finally, transvection of the above ($g_a(z_k)\wto 0$ and $g_{-a}(g_a(z_k))\wto 0$) leads to the conclusion that $a=0$ or $\norm{g_a(z_k)}\to 1$.  (Alternatively,  $g_a$ is norm continuous and fixes the boundary of the ball so that $\norm{z_k}\to 1$ implies $\norm{g_a(z_k)}\to 1$.)
\end{Proof}

   The following example shows that the conclusion in
   Proposition~\ref{prop:zero} above cannot be strengthened to ``$a=0$'' (or equivalently, that $g_a$ is linear).
\begin{Example}\label{ex:helpful}
  \begin{enumerate}
  \item Choose a basis of maximal tripotents of $J$ which are fixed by
    $j$.  (See Remark \ref{remark:basis}.)  Fix $t\in(0,1)$, $e=e_1$, $a=te$ and let
    $z_n=i\frac{n}{n+1}e_n$.  Then $z_n\wto 0$ and
    $\ip{z_n}{jz_n}\to -1$.  From \ref{cor:gmax},
    \[ g_a(z_k)= \left(t+
        (1-t^2)\frac{t\ip{z_k}{jz_k}}{1+t^2\ip{z_k}{jz_k}}\right)e +
      \frac{1-t^2}{1+t^2\ip{z_k}{jz_k}} i\frac k{k+1}e_k\] which also
    converges weakly to zero.

   Here, $g_a=g_{te}$ has $e$, and not $0$, as a fixed point.  This suggests that, generally, the fixed point of an automorphism (if it exists) 
    cannot be identified by using the above procedure to produce a weak fixed point and showing that the weak fixed point is indeed a fixed point.  This is verified below.
    
    \item With $J$ and $a=te$ as above, let $a_n = \frac 1n\to 0$ and $b_n=
      \sqrt{1-2ta_n+a_n^2}\to 1$.  Further let $\alpha_n= \frac 1{1+2ta_n}
      <1$ and $z_n= a_ne + ib_n e_n$.  The reader may verify that $\alpha_n g_a(z_n) = z_n$ and so
      $0=\wlim\ z_n $ is a weak fixed point of $g_a$, but is not a
      fixed point.
  \end{enumerate}
\end{Example}

Although this example shows that  a weak fixed point of $g$ may fail to be a fixed point, the following result gives conditions guaranteeing a weak fixed point is indeed fixed by $g$.

\begin{Theorem}\label{thm:first}
  Let $g=Tg_a$ be a biholomorphic mapping of the open unit ball $B$ of a spin factor $J$.
  Let $\xi$ be a weak fixed point of $g$ (that is, there are fixed points $z_k\in B$ of $\alpha_k g$ where
  $\alpha_k\nearrow 1$ and $\xi=\wlim z_k=\wlim g(z_k)\in \ol{B}$).   If any of the following conditions hold then $\xi$ is a fixed point of $g$.
  \begin{enumerate}[(1)]
  \item[(i)] $\norm \xi =1$, or
    \item[(ii)] $\xi=0$ and $a$ is not a non-zero scalar multiple of a
      maximal tripotent, or
    \item[(iii)] $\xi=0$ and $\norm{z_k}\not\to 1$, or
    \item[(iv)] each $z_k$ is rank one.
    \end{enumerate}

\end{Theorem}

\begin{Proof}
  First, suppose $\norm\xi=1$.  Without loss of generality, $z_k \wto
  \xi$ and $\norm{z_k}$ is convergent to $\ell$ say.  The Hahn-Banach
  theorem implies that $\ell=1$.  Thus $\frac{z_k}{\norm{z_k}} \wto
  \xi$ and so the
Kadec-Klee Property of the spin factor \cite[Corollary 3.6]{MR2176467} implies norm convergence of $z_k$
to $\xi$.  Norm continuity of the automorphism $g$ on $\overline{B}$ implies $\xi$ is a
fixed point of $g$ in $\ol B$.

In the case of (ii), we have $z_k=\alpha_k Tg_a(z_k)
\wto 0$ and thus $g_a(z_k) \wto 0$.  By Proposition~\ref{prop:zero}, $a=0$
and $g$ is linear, having a fixed point at $0$.  This proposition also covers case (iii).

In (iv),
  we have $z_k \wto \xi \in B$ and $Tg_a(z_k)=g(z_k)\wto \xi$.  By Lemma~\ref{lem:sfqi},
\[ z_k^{-a} = \frac 1{1+2\ip{z_k}a + \ip {z_k}{jz_k}\ip {ja}a}(z_k
+\ip{z_k}{jz_k}a).\]
Since each $z_k$ is rank one, $\ip{z_k}{jz_k}=0$ so that
\[z_k^{-a} = \frac 1{1+2\ip {z_k}a}z_k\]
which tends weakly to $\frac 1{1+2\ip \xi a} \xi=\xi^{-a}$.  It follows that
$g_{a}(z_k)=a+B(a,a)^{\half}(z_k^{-a}) $ has weak limit $g_a(\xi)$ and hence
$g(z_k)\wto g(\xi)$.
Thus $g(\xi)=\wlim g(z_k)=\wlim z_k=\xi$.



\end{Proof}

\begin{Remark}
    The proof of Theorem~\ref{thm:first} part (iv) above demonstrates that biholomorphic automorphisms are weakly continuous on rank one elements. 
\end{Remark}

We now determine conditions under which $g$ has a fixed point (that may not coincide with the weak fixed point which always exists).   Our proof is constructive and locates a fixed point on $\partial B.$

\begin{Theorem}\label{thm:orthogonal}
    Let $g=Tg_a$ have a weak fixed point at $0$ where $Tj=jT$ and the elements of $\{a, Ta, T^2a,\ldots\}$ are mutually orthogonal (that is, $\ip{T^ka}{T^ja}=0$ for $j\ne k$).  Then $g$ has a fixed point on $\partial B$.
\end{Theorem}

\begin{Proof}
    From Theorem~\ref{thm:first} part (ii), we can proceed to assume $a$ is a non-zero scalar multiple of a maximal tripotent, as otherwise, $0$ is a fixed point of $g$ and we are done.
    Therefore,  $a=\lambda e$, where $\lambda \in \Delta\setminus\{0\}$ and $e$ is a maximal tripotent.

    Without loss of generality, we assume that $je=e$.  Indeed,  Lemma~\ref{lem:maxtp} gives $je=\exp(i\theta) e$, for some $\theta\in \R$.  Now $e':=\exp(-i\theta/2) e$ satisfies $je'=e'$ and  (\ref{exp}) implies $g=Tg_{\lambda e}$ has a fixed point at $w$ (say) if, and only if, $Tg_{\lambda e'}$ has a fixed point (at $\exp(-i\theta/2)w$).  

    Similarly, there is no loss of generality in assuming $\lambda=t\in(0,1)$.  For this, let $\lambda=t\exp(i\phi)$ and note from (\ref{exp}) that $g=Tg_{\lambda e}$ has a fixed point at $y$ if and only if $Tg_{t e}$ has a fixed point at $\exp(-i\phi)y$.  
    
    Following these simplifications, the (Hilbert) orthogonality of $\{T^{n-1}a:n\in \N\}$ is equivalent to orthogonality of $\{T^{n-1}e : n\in \N\}$, and indeed orthogonality of $\mathcal B=\{T^n e: n\in\Z\}$, the elements of which are fixed by $j$.  In fact, $\ip{T^ke}{T^je}=\delta_j^k$ (Kronecker delta).  The set $\mathcal B$ may be extended to an orthonormal basis of $H$ (and $J$) fixed by $j$, but we will only need to consider elements in the (closed) span of $\mathcal B$.  

    From~Corollary~\ref{cor:gmax}, we know that 
    \[g_a(z)=g_{te}(z)= \left( t+ 
         \frac{t(1-t^2)\ip z{jz}}{1+2t\ip ze + t^2\ip z{jz}}\right)e 
         + 
         \left( \frac{1-t^2}{1+2t\ip ze + t^2\ip z{jz}}\right) z.
         \]
The scalar coefficients here of $e$ and $z$ will be denoted by $\lambda(z)$ and $\mu(z)$ respectively, giving,
\[g(z)=Tg_{te}(z)=  \lambda(z) Te 
         + 
         \mu(z) Tz.
         \]

Let $\mu_0:=\frac{1-t^2}{1+t^2}$ and $\lambda_0:=t(1+\mu_0)=\frac{2t}{1+t^2}$ and consider the element
\[ z_0:= \lambda_0 (Te + \mu_0 T^2 e + \mu_0^2 T^3e+\cdots),\]
which is well defined (since $|\mu_0|<1$ the sum converges in $J$).
We now proceed to prove that $z_0$ is our desired fixed point.

Clearly $jz_0=z_0$, as $\lambda_0$ and $\mu_0$ are real and $j$ fixes elements of $\mathcal B$.  Orthogonality of the terms allows us to calculate
\[ \ip {z_0}{jz_0}=\ip {z_0}{z_0} =\frac{|\lambda_0|^2}{1-|\mu_0|^2} = 1, \]
namely, $\norm {z_0} =1$.  Furthermore \[\lambda(z_0)= t+\frac{t(1-t^2).1}{1+2t.0+t^2} = \frac{2t}{1+t^2}=\lambda_0\] and 
\[\mu(z_0)=\frac{1-t^2}{1+2t.0+t^2}= \mu_0.\]
Finally,
\begin{align*}
     g(z_0)&=Tg_{te}(z_0)\\
  &= \lambda(z_0)Te+\mu(z_0)Tz_0 \\ 
  &= \lambda_0Te +\mu_0 T\left(\lambda_0 Te + \lambda_0 \mu_0 T^2 e + \lambda_0\mu_0^2 T^3e+\cdots\right)\\
  &= z_0
\end{align*}
and we have found the required fixed point (on $\partial B$).
\end{Proof}

We now strengthen {Theorem}~\ref{thm:orthogonal} significantly by replacing the orthogonality condition there with a much weaker, albeit more technical, condition.

\begin{Theorem}\label{thm:sliver_condition}
    Let $g=Tg_{te}$ have a weak fixed point at $0$ where $Tj=jT$.  Let $a_k:=\ip {T^ke}{e}$ and define $f:\Delta\to\C$ by $f(x)=\sum_{k=1}^{\infty}a_kx^k$.  
    If $$\limsup_{u\to 1^-} f(u)>-\half$$ then $g$ has a fixed point on $\partial B$.
\end{Theorem}

\begin{Remark}
    Note that $a_k\in\R$ since $\ol{a_k}=\ip e{T^ke}=\ip{j(T^ke)}{je}=\ip {T^ke}e=a_k$.  Also $|a_k|^2\le \ip{T^ke}{T^ke}\ip ee =1$ so $a_k\in[-1,1]$ and $f$ converges on $\Delta$.  If $\{e, Te, T^2e,\ldots\}$ is an orthogonal set then $f$ is identically zero and we recover Theorem~\ref{thm:orthogonal}.  Here, we can deduce the existence of a fixed point if, for example, $\ip{T^ke}e\ge 0$ for all $k$.
    \end{Remark}

\begin{Proof}(of Theorem~\ref{thm:sliver_condition}) 

\noindent The proof follows a roughly similar argument to that of Theorem~\ref{thm:orthogonal}.
For $u\in (0,1)$, let 
\begin{equation}z(u)=t(1+u)(Te+uT^2e+u^2T^3e+\cdots) \in J.\label{eq:z}\end{equation}  Evidently, $jz=z$ as $u$ is real and $j(T^ke)=T^ke$ for each $k$.  We calculate $\ip zz$, taking advantage of absolute convergence to reorder terms and the fact that $\ip{T^ke}{T^je}= \ip{T^{|j-k|}e}e=\ip e{T^{|j-k|}e}=a_{|j-k|}$, while $\ip ee=1$.  From\eqref{eq:z} we have
\begin{align*}
    \frac 1{t^2(1+u)^2}\ip zz &= \ip{\summ k1\infty u^{k-1}T^ke}{\summ k1\infty u^{k-1}T^ke} \\
      &= \ip{Te}{Te} + \ip{Te}{uT^2e} + \ip{Te}{u^2T^3e} + \ip{Te}{u^3T^4e}+\cdots \\
      &\qquad + \ip{uT^2e}{Te} + \ip{uT^2e}{uT^2e} + \ip{uT^2e}{u^2T^3e} + \ip{uT^2e}{u^3T^4e}+\cdots \\
        &\qquad + \ip{u^2T^3e}{Te} + \ip{u^2T^3e}{uT^2e} + \ip{u^2T^3e}{u^2T^3e} + \ip{u^2T^3e}{u^3T^4e}+\cdots \\
        &\qquad + \cdots \\ &\qquad \vdots \\
    &=  1 + a_1 u + a_2u^2 + a_3u^3 + \cdots \\
    &\ + a_1u + u^2 + a_1u^3 + a_2u^4 + \cdots \\
    &\ + a_2u^2 + a_1u^3 + u^4 + a_1 u^5 + \cdots \\
    &\ + \\ &\ \vdots
\end{align*}
We sum here over the diagonals.  The main diagonal sum is $1+u^2+u^4+\cdots = \frac 1{1-u^2}$.  The first upper- and lower-diagonal sums are each equal to $a_1u(1+u^2+u^4+\cdots)= \frac{a_1u}{1-u^2}$.  The second upper- and lower-diagonal sums combined give us $\frac{2a_2u^2}{1-u^2}$.  Proceeding in this way, we find the sum of all terms is $\frac 1{1-u^2}(1+2a_1u+ 2a_2u^2+ 2a_3u^3+\cdots) =\frac 1{1-u^2}(1+2f(u))$.  Therefore we have
\[\ip zz = t^2\frac{1+u}{1-u}\left(1 + 2f(u)\right)\]
and $\ip zz=1$ if, and only if, $h(u)=\frac1{t^2}$, where $h(u):=\frac{(1+u)(1+2f(u))}{1-u}$ is continuous on $(0,1)$.  The limsup hypothesis on $f$ in the statement guarantees $\limsup_{u\to 1^-} h(u)=\infty$ and since $\lim_{u\to 0^+} h(u)=1$, the intermediate value theorem then guarantees the existence of  $u_0\in(0,1)$ such that $h(u_0)=\frac1{t^2}$ and then $z_0:=z(u_0)$  satisfies $\ip {z_0}{z_0}=\ip {jz_0}{z_0} = \norm {z_0} =1$.

In addition, \begin{equation}f(u_0)=\half \left(\frac{1-u_0-t^2(1+u_0)}{t^2(1+u_0)}\right)\end{equation}
and \[\ip {z_0}e=t(1+u_o)\ip{Te+u_0T^2e+u_0^2T^3e+\cdots}e= t\frac{1+u_0}{u_0}f(u_0).\]
Now we have 
\begin{align*}
    \mu(z_0)&=\frac{1-t^2}{1+2t\ip {z_0}e  +t^2\ip z{jz}} \\
      &= \frac{1-t^2}{1+2t^2\left(\frac{1+u_0}{u_0}\right)f(u_0)+t^2}\\
      &= \frac{1-t^2}{1+ \frac{1+u_0}{u_0}\left(\frac{1-u_0-t^2(1+u_0)}{1+u_0}\right)+t^2}
\\
    &= u_0.
\end{align*} 
Similarly $\lambda(z_0)=t(1+u_0)$.  This gives
\begin{align*}
g(z_0)=Tg_{te}(z_0)&= \lambda(z_0)Te + \mu(z_0)Tz_0 \\
    &= t(1+u_0)Te + u_0Tz_0 \\
        &= t(1+u_0)Te + u_0 T\left( t(1+u_0)(Te+u_0T^2e+u_0^2T^3e+\cdots)\right) \\
        &= t(1+u_0)(Te + u_0T^2e+ u_0^2T^3e+\cdots) \\
        &= z(u_0)=z_0
\end{align*}
and we have found the desired fixed point $z_0$ with $\norm{z_0}=1$.
\end{Proof}

\section{Target set of a fixed-point free holomorphic map}

In this section we return to holomorphic maps without fixed points in $B$, namely let $f:B\to B$ be a compact, fixed-point free, holomorphic map on the open unit ball $B$ of a finite rank \JBT. It is well known, even for $B=\Delta\times \Delta$, that the sequence of iterates $(f^n)_n$ of $f$ need not converge. (Recall  $f$ compact means $f(B)$ is relatively compact.)

The set of accumulation points, $\Gamma(f)$, of $\{f^n: n \in \N\}$ (for any topology finer than the
topology of pointwise convergence on $B$ \cite[Corollary 3.9]{Mackey_Mellon_topologies}) and their images in $\ol B$ were examined by the current authors in \cite{MackeyMellon_WH}. The aim there was to locate the target set $T(f)$ of $f$ given by 
$$T(f):=\bigcup_{g\in \Gamma(f)} g(B),$$ which is necessarily a subset of $\partial B$.  It is shown in \cite[Theorem 2.10]{MackeyMellon_WH} that there is a unique closed convex set $W(f)$ in $\partial B$, termed the Wolff hull of $f$,  such that for all $x\in T(f)$, $$W(f)\cap Ch(x) \neq \emptyset$$
where $Ch(x):=\overline{K_x}$, and $K_x$ is the holomorphic boundary component of $x$. In this way, the Wolff hull $W(f)$ is proximal to the target set $T(f)$ and this provides important  limitations as to where on $\partial B$ we must look to find the target set.

For spin factors we can go further and  explicitly determine the location of the target set $T(f)$. In this case, we show that $T(f)$ is essentially two-dimensional in character. 

\begin{Theorem}\label{Wolff hull}
Let $B$ be the open unit ball of a spin factor $J$ and let $f:B\to B$ be a compact, fixed-point free, holomorphic map.
    Then there exists a  minimal tripotent $e\in \partial B$ (unique up to conjugation on $J$) such that $$T(f) \subseteq \partial B \cap (\C e \times \C \j e).$$
\end{Theorem}

\begin{Proof}
    Let $\xi$ be the "Wolff point" of $f$ as found in \cite{MackeyMellon_WH}.  Spectral decomposition of a boundary point of $B$ allows us to write $\xi=e+ \lambda \j e$ for a minimal tripotent $e$ and scalar $\lambda \in \ol \Delta$.   From \cite[Theorem~2.10]{MackeyMellon_WH}, we have that for all $x\in T(f)$, $$\ol {K_x} \cap \ol{K_{\tilde e}} \ne 0,$$  where since the rank is $2$ either $\tilde e=e$,    or $\tilde e=e+\mu \j e$  or $\tilde e = \mu \j e$ for $|\mu|=1$.  As $K_x=K_d$ for some tripotent $d$, we then have the existence of $v\in \ol{K_d}\cap \ol{K_{\tilde e}}$.

    As $\ol{K_d}=\cup_{p\ge d} K_p$ (and similarly for $K_{\tilde e}$), we have $v\in K_p\cap K_q$ for some tripotents $p\ge d$ and $q\ge \tilde e$.  In fact, as boundary components form a partition of $\partial B$ and each contain a unique tripotent, we conclude $p=q$.  Thus, there is a tripotent $p$ with $p\ge \tilde e$ and $p\ge d$.

    If $\tilde e = e$  or $\tilde e = \mu\j e$ then $p\ge \tilde e$ implies $$p\in\{e, e+\rho \j e : |\rho|=1, \mu \j e, \mu \j e+ \lambda e: |\lambda|=1\}.$$ On the other hand, if $\tilde e = e+\mu \j e$ then $p\ge \tilde e$ forces $$p=\tilde e=e+\mu\j e.$$  Considering both cases above, the only possibilities now for $d\le p$ are $$d \in \{\rho e, \mu \j e, \rho e+\mu\j e\}$$ with $|\rho|=|\mu|=1$.  The corresponding boundary components $K_d$ are $$\{\rho e+ \Delta \j e\}, \{\mu \j e + \Delta e\}, \ \hbox{and}\  \rho e+\mu \j e.$$  As $x\in K_x=K_d$, we have $x \in \C e \times \C \j e $, giving the stated result.
\end{Proof}

We now rephrase the above. Let $e \in J$ be any minimal tripotent and let $P_e$ be the vector subspace given by $$P_e:=\C e \times \C (je).$$ Since $e$ and $je$ are triple orthogonal from Lemma~\ref{lem:orthmintp} (see also Remark ~\ref{remark:basis}), it follows that $P_e$ is a $JB^*$-subtriple of $J$ and the spin factor norm of $J$ restricted to $P_e$ is just the maximum norm there. In other words, $P_e$ is  a copy of $\C^2$ with the maximum norm and the  open unit ball $B_e$  of $P_e$ is  a copy of the bidisc, namely, 
$$B_e=\Delta e \times \Delta (je).$$ This yields the following.
\begin{Corollary}\label{bidisc}
    Let $B$ be the open unit ball of a spin factor and let $f: B\to B$ be a compact, fixed-point free, holomorphic map.
    There exists a unique bidisc $B_e$ in $B$ such that $$T(f) \subseteq \partial B_e.$$
\end{Corollary}

\bibliographystyle{acm}

\def\cprime{$'$}

\end{document}